\documentclass[a4paper,11pt,reqno,oneside]{article}

\usepackage{amsmath}
\usepackage{amsthm}
\usepackage{amssymb}
\usepackage{enumerate}
\usepackage{graphicx}

\newcommand{\Z}{{\mathbb Z}}
\newcommand{\N}{{\mathbb N}}
\newcommand{\Q}{{\mathbb Q}}

\newcommand{\Acal}{{\mathcal A}}
\newcommand{\C}{{\mathcal C}}

\newcommand{\iet}{\ensuremath{\mathrm{iet}}}
\newcommand{\coloneq}{\mathrel{\mathop:}=}

\newtheorem{lem}{Lemma}
\newtheorem{thm}[lem]{Theorem}
\newtheorem{prop}[lem]{Proposition}
\newtheorem{coro}[lem]{Corollary}

\theoremstyle{definition}

\newtheorem{de}{\emph{Definition}}

\begin{document}

%%%%%%%%%%%%%%%%%%%%%%%%%%%%%%%%%%%%%%%%%%%%%%%%%%%%%%%%%%%%%%%%%%%%%%%%%%%%%

\title{\bf On the number of factors\\ in codings of three interval exchange}

\author{P. Ambro\v{z}$^{(1)}$, A. Frid$^{(2)}$, Z. Mas\'akov\'a$^{(1)}$, E. Pelantov\'a$^{(1)}$\\[2mm]
{\small $^{(1)}$ Doppler Institute \& Department of Mathematics}\\
{\small FNSPE, Czech Technical University in Prague}\\
{\small Trojanova 13, 120 00 Praha 2, Czech Republic}\\[2mm]
{\small $^{(2)}$ Sobolev Institute of Mathematics}\\
{\small Siberian Branch of the Russian Academy of Sciences}\\
{\small 4 Acad. Koptyug avenue, 630090 Novosibirsk, Russia}
 }

\date{}

\maketitle

\begin{abstract}
We consider exchange of three intervals with permutation
$(3,2,1)$. The aim of this paper is to count the cardinality of
the set $3\iet(N)$ of all words of length $N$ which appear as factors in infinite words coding such transformations.  We use
the strong relation of 3iet words and words coding exchange
of two intervals, i.e., Sturmian words. The known asymptotic
formula $\# 2\iet(N)/N^3\sim\frac1{\pi^2}$ for the number of
Sturmian factors allows us to find bounds $\frac1{3\pi^2} + o(1)
\leq \# 3\iet(N)/N^4 \leq \frac2{\pi^2} + o(1)$.
\end{abstract}

%\begin{keyword}
%Interval exchange transformation, Sturmian words, morphism, incidence matrix
%\MSC 68R15 \sep 15A36
%\end{keyword}
%\end{frontmatter}

%%%%%%%%%%%%%%%%%%%%%%%%%%%%%%%%%%%%%%%%%%%%%%%%%%%%%%%%%%%%%%%%%%%%%%%%%%%%%%%%%%%%%%%
\section{Introduction}

The study of infinite words arising from an exchange of several
intervals was initiated by Rauzy, whereas dynamical systems
connected with exchange of intervals has been already studied by
Katok and Stepin~\cite{StepinKatok}. The most explored subfamily
of such words are Sturmian words corresponding to an exchange of
two intervals \cite{fogg}. Infinite words coding exchange of 3
intervals form the next interesting family.

The special case of Sturmian words is very well understood.
%Numerous equivalent definitions have been given,
%for example as
%balanced words~\cite{morse-hedlund-ajm-62}, using
%palindromes~\cite{Droubay}, return words~\cite{vuillon}, and
%others. Criterions for substitutivity~\cite{????}
%and substitution invariance~\cite{yasutomi} of Sturmian words are known. A
%nice survey of these and many other
For surveys on the properties of Sturmian words see~\cite{lothaire2,fogg}.

Much less is known about codings of $r$ intervals
for $r>2$. Combinatorial characterization of the language of such
words was given first for $r=3$
in~\cite{ferenczi-holton-zamboni-jam-89}, then for general $r$
in~\cite{belov}. Boshernitzan and
Caroll~\cite{boshernitzan_caroll} have found a sufficient
condition for substitutivity of an infinite word coding exchange
of $r$ intervals. The necessity of this condition in the special
case $r=3$ and reverse-order permutation was demonstrated
in~\cite{adamczewski-jtnb-14}. A criterion for substitution
invariance of words coding exchange of 3 intervals with
permutation (3,2,1) is given in~\cite{abmp-integers-8}
and~\cite{balazi_masakova_pelantova_2}.

Actually, infinite words coding exchange of 3 intervals with
permutation (3,2,1), here called 3iet words, are subject of this
paper. They stand out among interval exchange words by
the fact that, similarly to Sturmian words, they can be
geometrically represented by cut-and-project sequences. This was
first remarked in~\cite{langevin-pmh-23}, then developed
in~\cite{gmp-jtnb-15}. Even more narrow connection between 3iet
words and Sturmian words is displayed in~\cite{abmp-integers-8}.
This connection allowed us to estimate the number of factors of 3iet words using the following
enumeration formula for Sturmian words, first given by Lipatov \cite{lipatov}; see also \cite{mignosi-tcs-82,berstel-pocchiola-ijac-3}.

\begin{thm}{\cite{lipatov,mignosi-tcs-82,berstel-pocchiola-ijac-3}}\label{thm:sturmfactors}
The number $\#2\iet(N)$ of Sturmian factors of length $N$ is given by
\[
\#2\iet(N)=1+\sum_{k=1}^N (N+1-k)\varphi(k) =
\frac{1}{\pi^2}N^3+{\mathcal O}(N^2\log N)\,,
\]
\end{thm}

The aim of this paper is to show the following theorem.

\begin{thm}\label{thm:hlavni}
The asymptotic growth of the number of 3iet factors of length $N$
is given by
\[
\frac{1}{3\pi^2}N^4 + {\mathcal O}(N^{3+\delta}) \leq
\#3\iet(N)  \leq \frac{2}{\pi^2}N^4+ {\mathcal
O}(N^{3+\delta})\,,
\]
for arbitrary $\delta>0$.
\end{thm}

The proof of Theorem~\ref{thm:hlavni} is divided into two parts,
the upper bound is given in Theorem~\ref{thm:horniodhad} in Section \ref{upper}, the
lower bound in Theorem~\ref{thm:dolniodhad} in Section \ref{lower}.
Lemmas on Sturmian words used in these sections may be of
independent interest.
Sections~\ref{totient} and \ref{sec:sturmian} describing properties of the Euler
function and of Sturmian factors are auxiliary.

Precise values of $\#3\iet(N)$ for $N\leq 10$ are found in Section
\ref{last}.
%%%%%%%%%%%%%%%%%%%%%%%%%%%%%%%%%%%%%%%%%%%%%%%%%%%%%%%%%%%%%%%%%%%%%%%%%%%%%%%%%%%%%%%
\section{Definitions of interval exchange words}\label{sec:preli}

Let $\Acal$ be a finite alphabet. A concatenation $w=w_0w_1\cdots
w_{n-1}$ of letters $w_i\in\Acal$ is called a word of length $|w|=n$.
The number of occurrences of a letter $a$ in the word $w$ is
denoted by $|w|_a$. The set of all finite words over $\Acal$ together
with the empty word is denoted by $\Acal^*$. We can also consider
one- or bidirectional infinite words $u=u_0u_1\cdots \in\Acal^\N$,
resp. $u=\cdots u_{-2}u_{-1}u_0u_1u_2\cdots \in\Acal^\Z$. A finite
word $w$ is a factor of an infinite word $u=(u_n)$ if
$w=u_iu_{i+1}\cdots u_{i+n-1}$ for some $i$. The factors of an
infinite word $u$ form the language of $u$, which is denoted by
${\mathcal L}(u)$. We denote ${\mathcal L}_n(u)={\mathcal L}(u)
\cap \Acal^n$. The factor complexity of $u$ is the function
$\C_u:\N\mapsto\N$ given by $\C_u(n)=\# {\mathcal L}_n(u)$.

Sturmian words are aperiodic words with the smallest possible
complexity, namely $\C_u(n)=n+1$ for all $n\in\N$. It is
known~\cite{morse-hedlund-ajm-62}
that Sturmian words on $\{0,1\}$ are exactly the lower and upper mechanical
words, which means that each Sturmian word $u=(u_n)_{n\in\Z}\in\{0,1\}^\Z$ is defined by
\begin{align*}
u_n & = \lfloor (n+1)\alpha+\beta\rfloor - \lfloor n\alpha+\beta\rfloor
\ \text{ or } \\
u_n &= \lceil (n+1)\alpha+\beta\rceil - \lceil n\alpha+\beta\rceil\,,
\end{align*}
where $\alpha\in(0,1)$ is irrational and $\beta\in[0,1)$. Another
representation of Sturmian words is given by the two interval
exchange. Define $\varepsilon=1-\alpha\in(0,1)$, where $\alpha$ is the parameter used above. Consider the
intervals $I_0=[0,\varepsilon)$, $I_1=[\varepsilon,1)$, $I_0\cup
I_1=[0,1)$, and the corresponding transformation $T_\varepsilon:
[0,1)\mapsto[0,1)$ defined by
\begin{equation}\label{eq:2iet}
T_{\varepsilon}(x)=
\begin{cases}
x+1-\varepsilon &\text{for } x\in I_0\,,\\
x-\varepsilon   &\text{for } x\in I_1\,.
\end{cases}
\end{equation}
We can code the orbit of a point $x_0\in [0,1)$ by the infinite
word $u_{\varepsilon,x_0}=(u_n)_{n\in\Z}\in\{0,1\}^{\Z}$:
\begin{equation}\label{eq:2ietslovo}
u_n=
\begin{cases}
0 & \text{if } T_{\varepsilon}^n(x_0)\in I_0\,,\\
1 & \text{if } T_{\varepsilon}^n(x_0)\in I_1\,.
\end{cases}
\end{equation}
The word $u_{\varepsilon,x_0}$ is the upper mechanical word corresponding to
$\alpha=1-\varepsilon$ and $\beta=x_0$. Similarly, lower
mechanical words are obtained by coding exchange of intervals of
the type $(\cdot,\cdot]$. It is the representation of Sturmian
words by two interval exchange which will be the most used here.

Since the language of a Sturmian word does not depend on the
initial point $x_0$, but rather only on the parameter
$\varepsilon$,  we can denote the language of $u_{\varepsilon,x_0}$ by
${\mathcal L}(\varepsilon)$ and the set of factors of length $M$
of $u_{\varepsilon,x_0}$ by ${\mathcal L}_M(\varepsilon)$.

We study infinite words which generalize Sturmian words
in the way that they code exchange of three intervals. Let
$\varepsilon,\ell$ be real numbers satisfying
\begin{equation}\label{eq:007}
\varepsilon\in(0,1)\setminus\Q\,,\qquad
\max\{\varepsilon,1-\varepsilon\}<\ell<1\,.
\end{equation}
The mapping $T_{\varepsilon,\ell}(x):[0,\ell)\mapsto [0,\ell)$
defined by
\begin{equation}\label{eq:1}
T_{\varepsilon,\ell}(x)=
\begin{cases}
x+1-\varepsilon & \text{for } x\in I_A\coloneq[0,\ell-1+\varepsilon)\,,\\
x+1-2\varepsilon &\text{for } x\in I_B\coloneq[\ell-1+\varepsilon,\varepsilon)\,,\\
x-\varepsilon &\text{for } x\in I_C\coloneq[\varepsilon,\ell)\,,
\end{cases}
\end{equation}
is called exchange of three intervals\footnote{Usually one defines
exchange of three intervals of lengths $\alpha,\beta,\gamma>0$
satisfying $\alpha+\beta+\gamma=1$. Our choice of
parameters~\eqref{eq:007} implies relation
$\alpha+2\beta+\gamma=1$, where $\varepsilon=\alpha+\beta$ and
$\ell=\alpha+\beta+\gamma$.} with permutation $(3,2,1)$. The orbit
of a point $x_0\in[0,\ell)$ under the transformation
$T_{\varepsilon,\ell}$ of~\eqref{eq:1} can be coded by the infinite
word $u_{\varepsilon,\ell,x_0}=(u_n)_{n\in\Z}$ in the alphabet
$\{A,B,C\}$, where
\begin{equation}\label{eq:2}
u_n=
\begin{cases}
A & \text{if } T_{\varepsilon,\ell}^n(x_0)\in I_A\,,\\
B &\text{if } T_{\varepsilon,\ell}^n(x_0)\in I_B\,,\\
C &\text{if } T_{\varepsilon,\ell}^n(x_0)\in I_C\,.
\end{cases}
\end{equation}
A visual representation of a 3 interval exchange transformation
can be found, e.~g., in \cite{adamczewski-jtnb-14}.

Since $\varepsilon$ is irrational, the infinite word $u$ is
aperiodic. We call these words 3iet words. Again, the language of
a 3iet word does not depend on the initial point $x_0$, thus we
denote it by ${\mathcal L}(\varepsilon,\ell)$. The factor
complexity of 3iet words depends on the parameter $\ell$ in the
following way: $\C(n)=2n+1$, if $\ell\not\in\Z+\varepsilon\Z$,
otherwise it satisfies $\C(n)=n+\mathit{const}$ for sufficiently
large $n$, see~\cite{adamczewski-jtnb-14}. Words with the latter
complexity are called quasisturmian and described for example
in~\cite{cassaigne,chernyatev}.

%%%%%%%%%%%%%%%%%%%%%%%%%%%%%%%%%%%%%%%%%%%%%%%%%%%%%%%%%%%%%%%%%%%%%%%%%%%%%%%%%%%%%%%

%%%%%%%%%%%%%%%%%%%%%%%%%%%%%%%%%%%%%%%%%%%%%%%%%%%%%%%%%%%%%%%%%%%%%%%%%%%%%%%%%%%%%%%
\section{Asymptotic behaviour of the totient
function}\label{totient}

In determining the asymptotic behaviour of the number of 3iet
factors, we shall strongly use asymptotic properties of the Euler
totient function,
$$
\varphi(n)=\#\{k\in\N\mid k\leq n,\ k\perp n\} =
n\prod_{\stackrel{p|n}{\hbox{\tiny $p$ prime}}}
\Big(1-\frac1{p}\Big)\,,
$$
where $k\perp n$ means $\gcd(k,n)=1$. In our formulas, we will use
the Landau big ${\mathcal O}$ notation; instead of a function
$f(n)$, we write ${\mathcal O}(g(n))$ if there exists a constant
$K$ such that $|f(n)| \leq \; K |g(n)|$ for all $n\in\N$.

The proof of the first asymptotic formula for the totient function can be found
for example in~\cite{nathanson-gtm-195}:
\begin{equation}\label{e:euler1}
\sum_{k=1}^n \varphi(k)= \frac{3}{\pi^2}n^2 + \mathcal{O}(n\log
n)\,.
\end{equation}
\iffalse
From the latter, we can easily derive
\begin{align*}
\sum_{k=1}^n k\varphi(k) &= \sum_{k=1}^n \sum_{i=1}^k \varphi(k) =
\sum_{i=1}^n\sum_{k=i}^n\varphi(k) = \sum_{i=1}^n
\Big(\sum_{k=1}^n\varphi(k) -
\sum_{k=1}^{i-1}\varphi(k) \Big) = \\[1mm]
&= \sum_{i=1}^n\Big(\frac{3}{\pi^2}n^2 - \frac{3}{\pi^2}i^2 +
\mathcal{O}(n\log n)\Big) = \frac{3}{\pi^2}n^3 -
\frac{1}{\pi^2}n^3 + \mathcal{O}(n^2\log n) = \\[1mm]
&= \frac{2}{\pi^2}n^3 + \mathcal{O}(n^2\log n)\,.
\end{align*}
It is then straightforward to obtain the growth order of the
number of Sturmian factors of length $N$,
\fi
Using the first equality from Theorem \ref{thm:sturmfactors} and this formula, one can easily derive the following estimate found already in \cite{lipatov}:
\begin{equation}\label{e:euler3}
\#2\iet(N)=1+\sum_{k=1}^N (N+1-k)\varphi(k)= \frac{1}{\pi^2}N^3 +
 \mathcal{O}(N^2\log N)\,.
\end{equation}

Yet another formula involving the totient function will be useful.
It can be easily derived that for any $q$, numbers coprime to $q$
are in a certain sense uniformly distributed, namely
$$
\# \{p\in\N \mid
 p\perp q,\ p\leq i\} = \frac{\varphi(q)}{q}i +
 {\mathcal O}(2^{\omega(q)})\,,
$$
where $\omega(q)$ is the number of prime divisors of $q$. From the
result of Hardy and Wright~\cite{hardy-wright-itn} on the asymptotic
behavior of $\omega(q)$, it follows that
\begin{equation} \label{e:euler4}
 \#\{p\in\N \mid p\perp q,\ j\leq p\leq i\} =
 \frac{\varphi(q)}{q}(i-j) + {\mathcal O}(q^{\delta})\,,\qquad
 \hbox{for any $\delta>0$}.
\end{equation}

%%%%%%%%%%%%%%%%%%%%%%%%%%%%%%%%%%%%%%%%%%%%%%%%%%%%%%%%%%%%%%%%%%%%%%%%%%%%%%%%%%%%%%%
\section{ Sturmian factors and Farey numbers}\label{sec:sturmian}
\iffalse
Factors of Sturmian words are certainly the most studied subject
of combinatorics on words~\cite{lothaire2,fogg}.
Here we recall only those properties which are important for
our study. It has been shown already by Hedlund and
Morse~\cite{morse-hedlund-ajm-62} that any language ${\mathcal
L}(\varepsilon)$ of a Sturmian word is closed under reversal,
i.e., $w_1w_2\cdots w_n\in{\mathcal L}(\varepsilon)$ implies
$w_nw_{n-1}\cdots w_1\in{\mathcal L}(\varepsilon)$. They have also
shown that ${\mathcal L}(\varepsilon)$ has the balance property:
for any pair of factors $w,v\in{\mathcal L}(\varepsilon)\subset
\{0,1\}^*$ of the same length we have $\bigl||w|_0-|v|_0\bigr|\leq
1$.
\fi
%is anything from this paragraph needed later?

As we have already mentioned, for every fixed irrational
$\varepsilon\in(0,1)$ and fixed $M\in\N$, the number of factors of
length $M$ in the language ${\mathcal L}(\varepsilon)$ of a
Sturmian word is exactly $M+1$. Obviously, the set ${\mathcal
L}_M(\varepsilon)$ does not determine $\varepsilon$ uniquely: the
same set of $M+1$ factors appears for uncountably many irrational
$\varepsilon$ which form an interval. As shown
in~\cite{mignosi-tcs-82}, if for some irrational
$\varepsilon_1<\varepsilon_2$ it holds that ${\mathcal
L}_M(\varepsilon_1)\neq{\mathcal L}_M(\varepsilon_2)$, then there
exists a rational $\frac{p}{q}$ such that
$\varepsilon_1<\frac{p}{q}<\varepsilon_2$ and $q\leq M$. Therefore
the interval $(0,1)$ of admissible parameters $\varepsilon$ is
partitioned into small intervals determining classes of Sturmian
words with distinct ${\mathcal L}_M(\varepsilon)$. The partition
is given by the Farey fractions of order $M$, i.e., all the
reduced fractions in $[0,1]$ with denominator smaller than or
equal to $M$. It is useful to consider the Farey fractions as an
ordered list
$$
{\cal F}_M: \quad
f_0=\frac{0}{1}<f_1<\cdots<f_{r-1}<f_r=\frac{1}{1}\,.
$$
%The Farey fractions have interesting properties. It is not
%difficult to see that the number of elements of ${\cal F}_M$ is
%equal to $1+r=1+\sum_{k=1}^M\varphi(k)$.
%% One can also show that
%% two consecutive fractions $\frac{p}{q}$,
%% $\frac{\tilde{p}}{\tilde{q}}\in{\mathcal F}_M$ satisfy
%% $\tilde{p}q-p\tilde{q}=1$. Conversely, reduced fractions
%% $\frac{p}{q}$, $\frac{\tilde{p}}{\tilde{q}}\in[0,1]$ with
%% $q,\tilde{q}\leq M$ are consecutive Farey fractions, if
%% $q+\tilde{q}\geq M+1$. Otherwise, the fraction
%% $\frac{p+\tilde{p}}{q+\tilde{q}}\in{\mathcal F}_M$ and
%% $\frac{p}{q}<\frac{p+\tilde{p}}{q+\tilde{q}}<\frac{\tilde{p}}{\tilde{q}}$.

\begin{figure}[!ht]
\centering\includegraphics{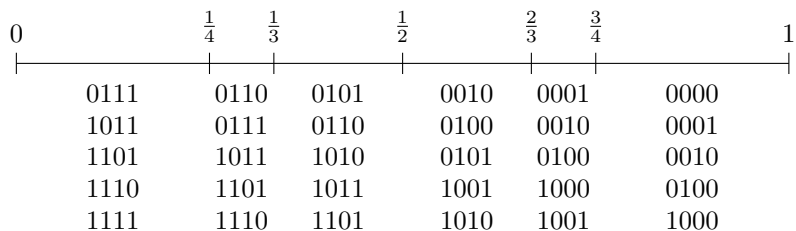} \caption{Classes of
Sturmian words of length 4.}
\label{f:sturmclasses}
\end{figure}

The partition $[0,1]=\cup_{i=0}^r[f_i,f_{i+1}]$, where ${\cal
F}_M=\{f_i\mid i=0,\dots,r\}$, is essential for describing
Sturmian factors of length $M$. We shall thus speak about classes
of Sturmian words of length $M$ determined by intervals
$(f_i,f_{i+1})$, $i=0,\dots,r-1$. Every class contains $M+1$
factors and of course, the same factor may belong to different
classes. For $M=4$, the division is illustrated in
Figure~\ref{f:sturmclasses} together with the corresponding sets
of factors of length $4$.

The possible numbers of letters 0 and 1 in factors corresponding to an interval
$(f_i,f_{i+1})$ can be easily determined as follows.

\begin{lem}\label{l:1}
Let $M\in\N$ be fixed, and let $f_i,f_{i+1}$ be two consecutive
Farey fractions of order $M$. Let $\varepsilon$ be an irrational
in $(f_i,f_{i+1})$. Then for every factor $w$ of length $M$ of a
Sturmian word with parameter $\varepsilon$ we have
$$
|w|_0 \in \big\{ \lfloor Mf_i\rfloor,\lfloor Mf_i\rfloor + 1
\big\}\,,\qquad |w|_1 \in \big\{ M-\lfloor Mf_i\rfloor, M-
\lfloor Mf_i\rfloor - 1 \big\}\,.
$$
\end{lem}

\begin{proof} Denote $\alpha=1-\varepsilon$. If $w=w_0w_1\cdots w_{M-1}$ is
a factor of a Sturmian word with parameter $\varepsilon$, there
must exist $\beta\in[0,1)$ such that
$$
w_i=\big\lfloor (i+1)\alpha+\beta \big\rfloor - \big\lfloor
i\alpha+\beta \big\rfloor \in\{0,1\}\,, \qquad\hbox{ for } i
=0,1,\dots,M-1\,.
$$
Therefore
$$
|w|_1=\sum_{i=0}^{M-1}w_i = \lfloor M\alpha+\beta \rfloor \in
\big\{\lfloor M\alpha \rfloor,\,\lfloor M\alpha \rfloor+1\big\}\,.
$$
The language ${\mathcal L}_M(\varepsilon)$ is the same for all
$\varepsilon\in(f_i, f_{i+1})$ and we have $\lfloor M\alpha
\rfloor = \big\lfloor M(1-\varepsilon) \big\rfloor = M-\lceil
M\varepsilon\rceil = M-\lfloor M\varepsilon\rfloor-1 = M-\lfloor
Mf_i\rfloor-1$. \end{proof}

Now consider words coding the
transformations $T_{f_i}$ defined as in \eqref{eq:2iet}. If $f_i=\frac{p}{q}$,
then such a word is purely
periodic and looks as $v^{(i)}v^{(i)}v^{(i)}\cdots =
\big(v^{(i)}\big)^\omega$, where $v^{(i)}$ is a primitive word of
length $q$. Recall that a finite word $v$ is said to be primitive
if it cannot be written in the form $v=ww\cdots w=w^k$ for any
integer $k\geq 2$. Therefore the language of the periodic word
coding the transformation $T_{f_i}$ has exactly $q$ factors of
length $M$ for every $M\geq q$. Let us denote the set of such
factors by $L_M^{(i)}$. As it is derived
in~\cite{mignosi-tcs-82},
$$
{\mathcal L}_M(\varepsilon)=L_M^{(i)}\cup L_M^{(i+1)}\,, \quad
\hbox{for all }\varepsilon\in(f_i,f_{i+1})\,.
$$
We have $\#{\mathcal L}_M(\varepsilon)=M+1$, $\#L_M^{(i)}=q_i$,
and $\#L_M^{(i+1)}=q_{i+1}$, where $q_i,q_{i+1}$ are denominators
of Farey fractions $f_i, f_{i+1}$, respectively. Therefore we have
\begin{equation}\label{eq:prunik}
\# \Big(L_M^{(i)}\cap L_M^{(i+1)}\Big) = q_i+q_{i+1} - M-1\,.
\end{equation}

\section{The upper bound}\label{upper}

The main tool which we use in this section is a link
between 3iet words and Sturmian words over the alphabet $\{0,1\}$
by means of morphisms
$\sigma_{01},\sigma_{10}:\{A,B,C\}^*\mapsto\{0,1\}^*$ defined by
\begin{equation}\label{eq:002}
\begin{split}
\sigma_{01}:&\quad A \mapsto 0\,,\ B  \mapsto 01\,,\  C \mapsto 1\,,\\
\sigma_{10}:&\quad A \mapsto 0\,,\ B  \mapsto 10\,,\  C \mapsto 1\,.
\end{split}
\end{equation}
In~\cite{abmp-integers-8}, the following statement is proved.

\begin{thm}\label{thm:1}
A ternary word $u$ is a 3iet word if and only if both
$\sigma_{01}(u)$ and $\sigma_{10}(u)$ are Sturmian words.
\end{thm}

Based on this strong relation, we define the notion of
$b$-amicability.

\begin{de}
Let $w^{(1)},w^{(2)}$ be finite words over the alphabet $\{0,1\}$.
We say that $w^{(1)}$, $w^{(2)}$ form a $b$-amicable pair, if
there exists a 3iet factor $w$ over $\{A,B,C\}$ with exactly $b$
letters $B$ such that $w^{(1)}=\sigma_{01}(w)$ and
$w^{(2)}=\sigma_{10}(w)$.
\end{de}

\begin{figure}[!ht]
\begin{equation*}
\begin{aligned}
w^{(1)} =\quad \\[1mm] w^{(2)} =\quad \\[1mm] w =\quad
\end{aligned}
\begin{gathered}\boxed{\:\begin{gathered}0 \\ 0\end{gathered}\:}\\A\end{gathered}\
\begin{gathered}\boxed{\:\begin{gathered}1 \\ 1\end{gathered}\:}\\C\end{gathered}\
\begin{gathered}\boxed{\:\begin{gathered}0 \\ 0\end{gathered}\:}\\A\end{gathered}\
\begin{gathered}\boxed{\:\begin{gathered}0\ 1 \\ 1\ 0\end{gathered}\:}\\B\end{gathered}\
\begin{gathered}\boxed{\:\begin{gathered}0 \\ 0\end{gathered}\:}\\A\end{gathered}\
\begin{gathered}\boxed{\:\begin{gathered}1 \\ 1\end{gathered}\:}\\C\end{gathered}
\end{equation*}
\caption{Finite words $w^{(1)}=0100101$ and $w^{(2)}=0101001$
forming a $1$-amicable pair. The corresponding ternary word is
$w=ACABAC$.} \label{fff}
\end{figure}

The notion is illustrated on Figure~\ref{fff}. Note that due to
Theorem~\ref{thm:1}, both words of a $b$-amicable pair are
Sturmian factors from the same language ${\mathcal
L}(\varepsilon)$ for some irrational $\varepsilon\in(0,1)$. A 3iet
factor $w$ of length $N$ with precisely $b$ letters $B$
corresponds to a $b$-amicable pair of Sturmian factors of length
$N+b$. Obviously, different $b$-amicable pairs give rise to
different 3iet factors. Therefore we can express the number of
3iet words as the number of corresponding $b$-amicable pairs.

\begin{prop}\label{p:zaklad}
Let $N\in\N$. Then
$$
\#{3\iet}(N) = \sum_{b=0}^N \#3\iet(N,b) = \sum_{b=0}^N \#\text{
$b$-amicable pairs of length $(N+b)$}\,,
$$
where $3\iet(N,b)$ denotes the set of 3iet factors of length $N$
with precisely $b$ letters $B$.
\end{prop}

Let us focus on $b$-amicable
pairs of Sturmian factors. First, it is obvious that words which form a $b$-amicable
pair must both contain at least $b$ letters 0 and $b$ letters 1.
As a consequence of Lemma~\ref{l:1}, we have the following
statement.

\begin{coro}\label{c}
Let $b,M\in\N$ and let $\varepsilon$ be an irrational in $(f_i,
f_{i+1})$, where $f_i,f_{i+1}$ are consecutive Farey fractions of
order $M$. If ${\mathcal L}_M(\varepsilon)$ contains a
$b$-amicable pair, then
\[
b\leq \min \Big\{\lfloor Mf_i\rfloor+1, M-\lfloor
Mf_i\rfloor\Big\}\,.
\]
\end{coro}

\begin{lem}\label{l:3}
Let $Z_{M,b}$ be the family of classes of Sturmian words of length
$M$ containing at least one $b$-amicable pair. Then for $M< 2b$ we
have $\# Z_{M,b}=0$, and otherwise,
\[
\# Z_{M,b}\: \leq \frac3{\pi^2} (M-2b) M  + {\mathcal
O}(M^{1+\delta})\,,\quad \hbox{for arbitrary } \delta>0\,.
%
%\leq \ \frac{M-2b+2}{M} \sum_{q=1}^M \varphi(q) \ +\ {\mathcal
%O}(M^{1+\delta})\,,\quad \hbox{for some } \delta>0\,.
\]
\end{lem}

\begin{proof} Obviously, words forming a $b$-amicable pair must be of length
at least $M \geq 2b$.
%In order to find the bound on $\# Z_{M,b}$,
%realize that factors forming a $b$-amicable pair must have at
%least $b$ letters 0 and at least $b$ letters 1.
Consider an interval $[f_i,f_{i+1})$ such that the language
${\mathcal L}_M(\varepsilon)$ for $\varepsilon\in(f_i,f_{i+1})$
contains an amicable pair and $f_i=\frac{p}{q}$, $p\perp q$.
Corollary~\ref{c} then implies that $b$ must satisfy
\[
b\leq M\frac{p}{q}+1\qquad\hbox{and}\qquad b\leq M-M\frac{p}{q}+1\,,
\]
providing bounds on $p$, namely,
\[
q\frac{b-1}{M}\leq p \leq q\left (1-\frac{b-1}{M}\right )\,.
\]
Thus, we are looking for the number
\[
\sum_{q=1}^M \#\Bigl\{p\in\N \,\Big|\, p\perp q \hbox{ and }
q\tfrac{b-1}{M}\leq p \leq q\bigl(1-\tfrac{b-1}{M}\bigr)\Bigr\}
\]
With the use of~\eqref{e:euler4} we can write
\begin{align*}
\# Z_{M,b} & \leq \sum_{q=1}^{M} \frac{\varphi(q)}{q} q
\Big(1-\frac{2(b-1)}{M}\Big) + M{\mathcal
O}\big(M^{\delta}\big) =\\
& =  \frac{M-2b+2}{M} \sum_{q=1}^M \varphi(q) + {\mathcal
O}(M^{1+\delta})\, \mbox{~for all~} \delta>0.
\end{align*}
Together with the asymptotic behavior of $\sum_{q=1}^M \varphi(q)$
given in~\eqref{e:euler1}, this implies the statement of the lemma.
\end{proof}

\begin{lem}\label{l:2}
Let $\varepsilon\in(0,1)$ be irrational. Let $b,M\in\N$,
$2b\leq M$. Then the number of distinct $b$-amicable pairs in class
${\mathcal L}_M(\varepsilon)$ is at most $M-b$.
\end{lem}

\begin{proof}
 The language ${\mathcal L}_M(\varepsilon)$ contains $M+1$ factors of length
 $M$. Choose such  a factor $w$ and consider the set  of points $x\in[0,1)$ such
 that the orbit
$\{T^j_{\varepsilon}{(x)}\}$, $j=0,1,\dots,M-1$, is coded by $w$.
It can be easily shown that this set is a subinterval of $[0,1)$
of the type $[\cdot,\cdot)$. Thus $[0,1)$ is divided into $M+1$
 disjoint subintervals; the division is done by $M$ division points
 $T^{-j}_{\varepsilon}(\varepsilon)$, $j=0,1,\dots,M-1$. Let us
 denote these intervals $J_i$, $i=1,2,\dots,M+1$, ordered so that
$x<y$ for any $x\in J_i$ and any $y\in J_{i+1}$. The corresponding
words of length $M$ are denoted by
 $w^{(i)}$, $i=1,2,\dots,M+1$. Note that $T^{j}(J_i)$ is an
 interval for all $j=0,1,\dots,M-1$.

Let us first find an upper bound for the number of 1-amicable
pairs of Sturmian factors of length $M$. The following procedure is
illustrated on Figure~\ref{f:sturm_iter}. Take two consecutive
intervals, $J_i$, $J_{i+1}$ for some $i=1,\dots,M$. It means that
$J_i\cup J_{i+1}$ is an interval which contains in its interior
exactly one of the division points, say
$T^{-k}_\varepsilon(\varepsilon)$. If $k=0$, then the word
$w^{(i)}$ starts with 0 and the word $w^{(i+1)}$ starts with 1.
More generally, realize that words $w^{(i)}$ and $w^{(i+1)}$ have
a common prefix $v$ of length $k$ and the set
$T^{k}_\varepsilon(J_i\cup J_{i+1})$ is still an interval. Since
$T^{-k}_\varepsilon(\varepsilon)\in (J_i\cup J_{i+1})^\circ$, we
have $\varepsilon\in T^{k}_\varepsilon\big((J_i\cup
J_{i+1})^\circ\big)$, and therefore the word $w^{(i)}$ has prefix
$v0$ and the word $w^{(i+1)}$ has prefix $v1$.

If $k=M-1$, then $w^{(i)}=v0$ and $w^{(i+1)}=v1$ and they are not
1-amicable. If, by contrary, $k<M-1$, then
$T^{k+1}(J_i)=T^k(J_i)+1-\varepsilon$ is an interval of the form
$[\cdot,1)$, and $T^{k+1}(J_{i+1})=T^k(J_{i+1})-\varepsilon$ is an
interval of the form $[0,\cdot)$. We further have
\begin{align*}
T^{k+2}(J_i) & = T^{k+1}(J_i)-\varepsilon = T^k(J_i)+1-2\varepsilon\,,\\
T^{k+2}(J_{i+1}) & = T^{k+1}(J_{i+1})+1-\varepsilon = T^k(J_{i+1})+1-2\varepsilon\,.
\end{align*}
and therefore the set $T^{j}(J_i\cup J_{i+1})$ is again an
interval for all $j$, $k+2\leq j\leq M-1$, and thus the words
$w^{(i)}$ and $w^{(i+1)}$ have a common suffix $v'$ of length
$M-k-2$. Therefore we have $w^{(i)}=v01v'$ and $w^{(i+1)}=v10v'$
and the words $w^{(i)}$ and $w^{(i+1)}$ form a 1-amicable pair.
Together, we have a $1$-amicable pair of words $w^{(i)}$ and
$w^{(i+1)}$ for $i=1,\dots,M$, $i\neq k$. Their number is
therefore $M-1$.

\begin{figure}[!ht]
\centering\includegraphics{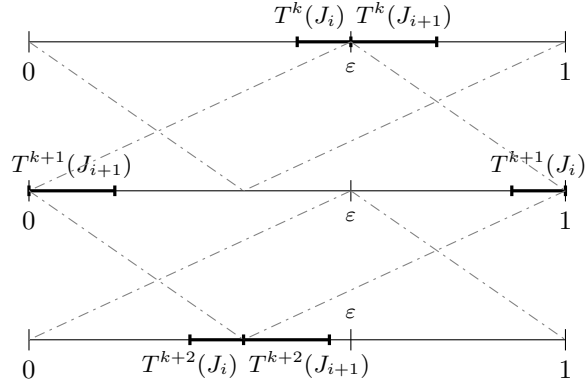} \caption{Iterations
of consecutive intervals $J_i$, $J_{i+1}$ show that the
corresponding words $w^{(i)}$, $w^{(i+1)}$ form a 1-amicable
pair.}\label{f:sturm_iter}
\end{figure}

The previous considerations imply that the words $w^{(i)}$,
$i=1,\dots,M+1$, are ordered lexicographically, where consecutive
words differ only by interchange of $01\leftrightarrow10$ at one
place. The only exception is the pair of words $w^{(i)}$ and
$w^{(i+1)}$ whose intervals $J_i$, $J_{i+1}$ are separated by the
iteration $T^{-M+1}_\varepsilon(\varepsilon)$. These words are of
the form $w^{(i)}=v0$ and $w^{(i+1)}=v1$, which corresponds to
`half' of the interchange $01\leftrightarrow10$.

From that, it is clear that $b$-amicable pairs can only be the
pairs of words $w^{(i)}$, $w^{(i+b)}$, $i=1,\dots,M+1-b$, with
exceptions, namely those indices $i$, such that the interval
$J_i\cup\dots\cup J_{i+b}$ contains in its interior the iteration
$T^{-M+1}_\varepsilon(\varepsilon)$. However, since we do not know
the position of this point, it is difficult to determine the
number of these exceptions, but it is at least one. Therefore the
number of $b$-amicable pairs is at most $M-b$. \end{proof}

Now we are ready to prove the upper bound from Theorem
\ref{thm:hlavni}.
\begin{thm}\label{thm:horniodhad}
The number of 3iet factors of length $N$ satisfies
$$
\# 3\iet(N)\leq \frac{2}{\pi^2} N^4 + {\mathcal O}(N^{3+\delta})\,
$$
for arbitrary $\delta>0$.
\end{thm}

\begin{proof}
 %For counting the 3iet factors we shall use
 Combining Proposition~\ref{p:zaklad}
% . Thus we are interested in counting
 %all Sturmian $b$-amicable pairs of length $N\!+\!b$. It is
 %difficult to determine this number exactly, since some classes of
 %Sturmian words of length $M$ do not contain any $b$-amicable
 %pair, and some $b$-amicable pairs are found in more than one
 %class. However, we can get an upper estimate, summing the upper
 %bound $M-b$ from
 and Lemma~\ref{l:2},
 %over all classes that have
 %at least one $b$-amicable pair of length $M$. In Section~\ref{sec:sturmian},
 %we have denoted the family of such classes by $Z_{M,b}$ and
 %derived an upper bound on their number in Lemma~\ref{l:3}.
we see that
%Realize that factors forming a $b$-amicable pair must have at
%least $b$ letters 0 and at least $b$ letters 1. Thus $Z_{M,b}$ is
%bounded from above by the XXXXXXXXXXXXXXXXXXXXXX
%  (their number is $\Z_{M,b}$ as
% given in Lemma~\ref{l:3}).
$$
\begin{aligned}
\# {3\iet}(N) &= \sum_{b=0}^N \# 3\iet(N,b) =
\sum_{b=0}^N \# \hbox{ $b$-amicable pairs of length $(N+b)$} \leq\\
&\leq \sum_{b=0}^N\  \sum_{{\mathcal K}\in Z_{N+b,b}}
\hspace*{-4mm}\# \hbox{ $b$-amicable pairs of length $(N+b)$ in the class ${\mathcal K}$}\leq\\
&\leq \sum_{b=0}^N N \# Z_{N+b,b}\,.
\end{aligned}
$$
Putting in the estimate on $\#Z_{N+b,b}$ from Lemma~\ref{l:3},
we obtain
\begin{align*}
\# {3\iet}(N) &\leq N \sum_{b=0}^N \frac{3}{\pi^2}(N-b)(N+b) +
N \sum_{b=0}^N {\mathcal O}\big((N+b)^{1+\delta}\big)=\\[2mm]
&= \frac{3}{\pi^2} N \Big( N^2(N+1) - \frac16 N(N+1)(2N+1)\Big)
+  {\mathcal O}(N^{3+\delta}) = \\[2mm]
&= \frac2{\pi^2}N^4 + {\mathcal O}(N^{3+\delta})\,.\qedhere
\end{align*}
 \end{proof}

\section{The lower bound}\label{lower}

Consider the following geometric
representation of Sturmian and 3iet words by cut-and-project
sequences. Given a strip in the Euclidean plane parallel to the
(irrationally oriented) straight line $y=\varepsilon x$, take all
points of the lattice $\Z^2$ and project them orthogonally to the
straight line. It is known~\cite{langevin-pmh-23,gmp-jtnb-15} that in such a way, one
obtains a sequence of points with two or three distances between
them. If the distances are two, say $\Delta_1$, $\Delta_2$, we can
code them by two letters, and the resulting word is a Sturmian word.
If the distances are three, they are of the form $\Delta_1$,
$\Delta_2$, $\Delta_1+\Delta_2$. Coding these distances by three
letters, we obtain a 3iet word. Moreover,
Sturmian word is obtained only for a discrete set of values of the
width of the considered strip. For details about this construction
we refer to~\cite{gmp-jtnb-15}.

On the other hand, any Sturmian word $u_{\varepsilon,x_0}$ or a 3iet
word $u_{\varepsilon,\ell,x_0}$ can be represented in this way, the
parameter $\varepsilon$ being the slope of the straight line and
$\ell$ corresponding to the width of the strip. If without loss of
generality we take $\varepsilon\in(0,1)$, the widths of the
strips giving rise to Sturmian words form a sequence
$$
\cdots < \max\{\varepsilon,1-\varepsilon\}  < 1  < 1+\varepsilon  < \cdots
$$
%Why a sequence, not a set? What about $2-\varepsilon$?
All parameters $\ell$ in between these values yield 3iet words.

\begin{figure}[!ht]
\centering
\includegraphics{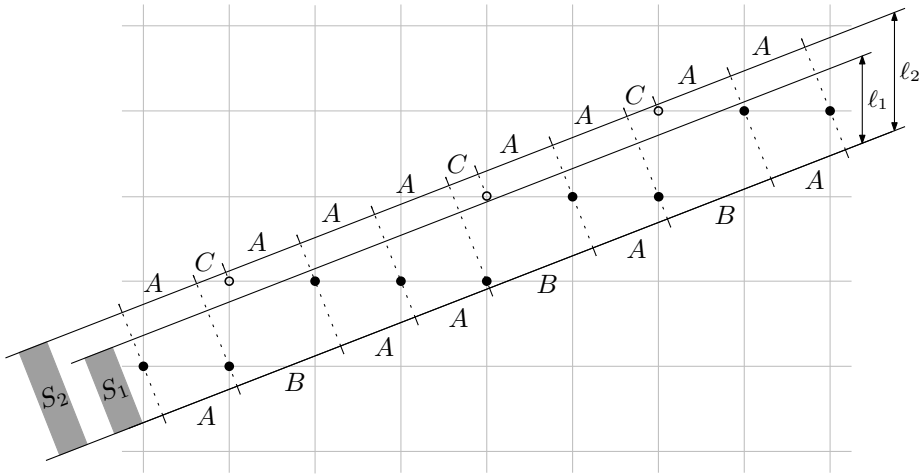}
\caption{Geometric representation of Sturmian and 3iet words by
cut-and-project sequences.} \label{f:projekce}
\end{figure}

This geometric representation is illustrated in
Figure~\ref{f:projekce}, where $\ell_1<\ell_2$ are two consecutive
discrete values for which the cut-and-project scheme yields a
Sturmian word. Projections of points in the strip $S_1$ with
parameter $\ell_1$ correspond to a Sturmian word $u^{(1)}$ over
the alphabet $\{A,B\}$, where $A$ stands for the projection of the
horizontal side and $B$ of the diagonal of a unit square. The strip
$S_2$ corresponding to the value $\ell_2$ gives a Sturmian word
$u^{(2)}$ over $\{A,C\}$, where again, $A$ is the projection of
the horizontal side and $C$ of the vertical side of the unit
square. Note that the projection of a diagonal is the sum of
projections of the vertical and horizontal sides. This corresponds
to the fact that the infinite word $u^{(2)}$ arises from $u^{(1)}$
replacing every $B$ by $AC$.

\begin{figure}[!ht]
\centering
\parbox{0.2\textwidth}{\centering\includegraphics[width=0.2\textwidth]{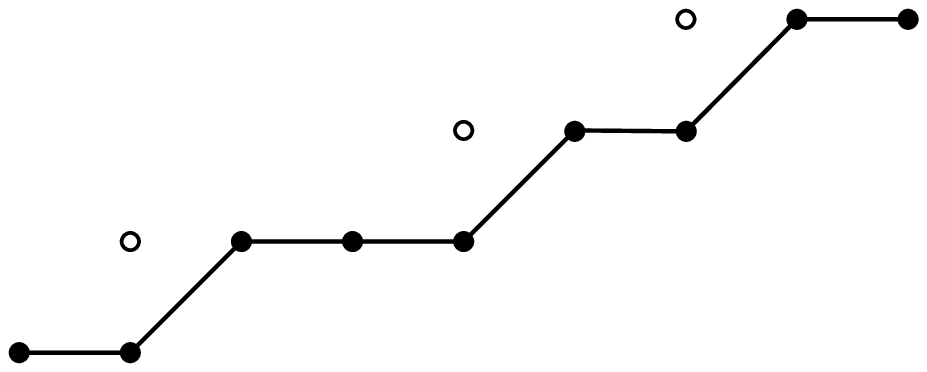}\par{\small $ABAABABA$}}
$\ \Rightarrow\ $
\parbox{0.2\textwidth}{\centering\includegraphics[width=0.2\textwidth]{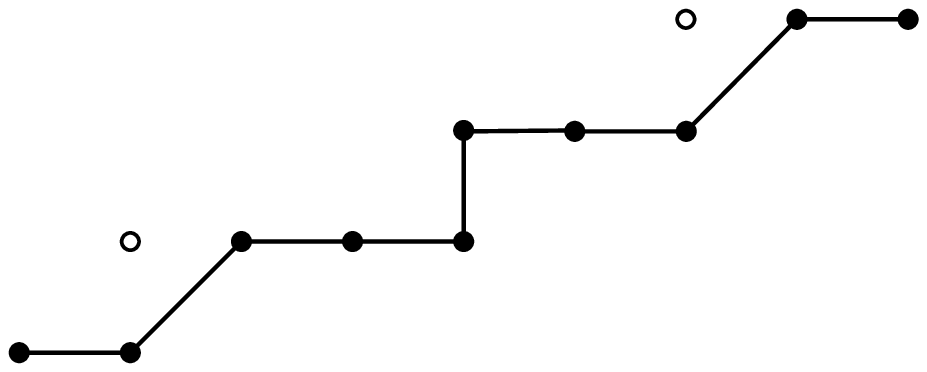}\par{\small $ABAACAABA$}}
$\ \Rightarrow\ $
\parbox{0.2\textwidth}{\centering\includegraphics[width=0.2\textwidth]{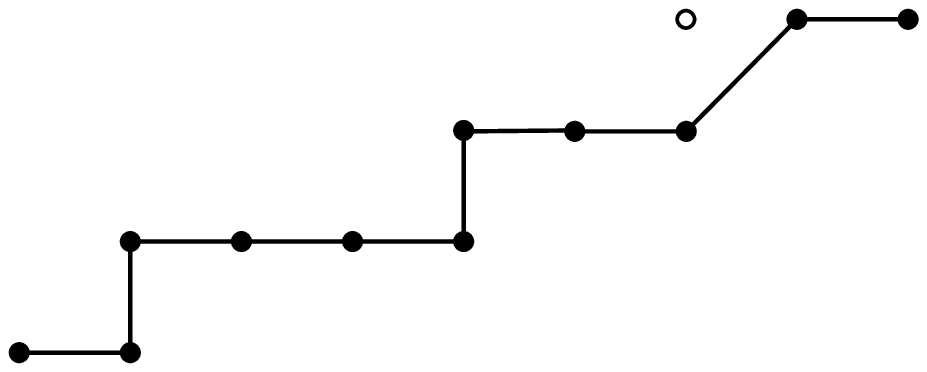}\par{\small $ACAAACAABA$}}
$\ \Rightarrow\ $
\parbox{0.2\textwidth}{\centering\includegraphics[width=0.2\textwidth]{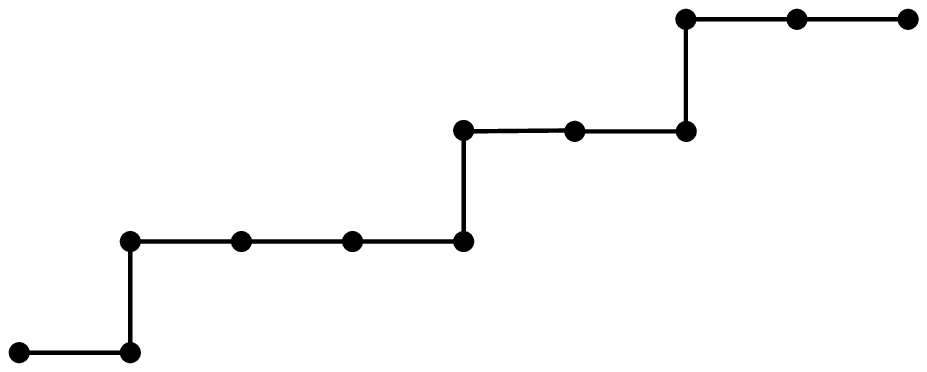}\par{\small $ACAAACAACAA$}}
\vspace*{2mm} \caption{Construction of 3iet factors with given
number of letters $C$ by substituting $AC$ for suitable letters
$B$ in a Sturmian factor.}
 \label{f:3ietpostupne}
\end{figure}

Any strip $S$ such that $S_1\subset S\subset S_2$ gives rise to a
3iet word over the alphabet $\{A,B,C\}$ where only some letters
$B$ have been replaced by $AC$. Enlarging the strip $S$ from $S_1$
to $S_2$ causes more and more lattice points to fall into $S$ and
thus to split more and more distances $B$ into $CA$. This procedure
is illustrated in Figure~\ref{f:3ietpostupne} on a finite segment
taken from Figure~\ref{f:projekce}.

Let now $w^{(1)}\in\{A,B\}^*$ be a Sturmian factor of length $M$
with at least $b$ letters $B$, $b\leq M$. By a suitable choice of
the enlarged strip $S$, we obtain from $w^{(1)}$ a
3iet factor $w$ of length $M+b$ containing exactly $b$ letters
$C$. Such a word $w$ satisfies $|w|_A\geq |w|_C$, where the equality
holds only if $w=(CA)^M$. Since the role of letters $A$ and $C$ is
symmetric over the set of all 3iet factors, the interchange
$A\leftrightarrow C$ in the factor $w$ gives again a 3iet factor.
Denoting by $H_{M,b}$ the number of Sturmian factors of length $M$
with at least $b$ letters $B$, we can summarize the above
considerations by saying that $\#3\iet(N)\geq 2\sum_{b=0}^N
H_{N-b,b}$. Moreover, realize that a Sturmian factor
with at least $b$ letters B is of length at least $b$, so the
sum is taken only over $b$ such that $N-b\geq b$, i.e.,
$b\leq \lfloor N/2\rfloor$. This gives the following

\begin{prop}\label{p:zaklad2}
Let $N\in\N$. Then
$$
\# 3\iet(N) \geq 2\sum_{b=0}^{\lfloor N/2\rfloor} \# \hbox{
Sturmian factors of length $(N-b)$ with at least $b$ letters
B}\,.
$$
\end{prop}

%Both propositions~\ref{p:zaklad} and~\ref{p:zaklad2} transform the question about the number of 3iet factors to a
%study of Sturmian factors. Thus, a detailed study of factors of
%Sturmian words will be needed (Section~\ref{sec:sturmian}). Since
%the overall number Sturmian factors is given using the totient
%function $\varphi(N)$,
%the following section is devoted to its asymptotic behavior.

%%%%%%%%%%%%%%%%%%%%%%%%%%%%%%%%%%%%%%%%%%%%%%%%%%%%%%%%%%%%%%%%%%%%%%%%%%%%%%55

To estimate the number in the right part of the equation,
we use the relation \eqref{eq:prunik}. Let us
mention that an asymptotic formula for the number of Sturmian
factors $w$ with given $|w|$ and $|w|_1$ was derived
in~\cite{bdjr-dmtca}.

\begin{lem}\label{l:4}
The number $H_{M,b}$ of Sturmian factors of length $M$ with at
least $b$ letters 1 satisfies
$$
H_{M,b} \geq \frac{1}{\pi^2} (M-b)M^2 + {\mathcal
O}(M^{2+\delta})\,,\qquad\hbox{for any }\delta>0\,.
$$
\end{lem}

\begin{proof}
Let $f_i$ be the $i$-th Farey
fraction $f_i=\frac{p_i}{q_i}$. Consider the set $L^{(i)}_M$. By Lemma~\ref{l:1}, if $b\leq
M-Mf_i-1$, then every factor in $L^{(i)}_M$ has at least $b$
letters 1.
The inequality $M-Mf_i-1\geq b$ can be rewritten as
\begin{equation}\label{eq:012}
p_i\leq \frac{M-b-1}{M}q_i\,.
\end{equation}
By~\eqref{e:euler4}, the number of Farey fractions with the same
denominator $q=q_i\leq M$ satisfying~\eqref{eq:012} is equal to
\begin{equation}\label{eq:013}
\frac{M-b-1}{M}\varphi(q) + {\mathcal O}(q^\delta)\,,\qquad \hbox{
for any } \delta>0\,.
\end{equation}

Since the values of Farey fractions increase with $i$,
there exists a maximal index, say $s$, such that every factor in
the set $\bigcup_{i=0}^s L^{(i)}_M$ has sufficient number of
letters 1. Therefore
\[
H_{M,b} \geq \# \bigcup_{i=0}^s L^{(i)}_M\,.
\]
In order to determine the cardinality of this union, we write it
as a disjoint union, namely
\[
\bigcup_{i=0}^s L^{(i)}_M = L^{(0)}_M\ \cup\ \bigcup_{i=1}^{s}
 \Big(L^{(i)}_M \setminus L^{(i-1)}_M \Big)\,.
\]
Since $\# L^{(i)}_M=q_i$ and  $\# L^{(i-1)}_M=q_{i-1}$ and
equality~\eqref{eq:prunik} holds, we obtain
\[
\# \Big(L^{(i)}_M \setminus L^{(i-1)}_M \Big) = M+1-q_{i-1}\,,
\]
which implies
\[
H_{M,b} \geq 1+ \sum_{i=1}^{s} (M+1-q_{i-1}) = \sum_{i=0}^{s} (M+1-q_{i}) - M + q_s\,.
\]
We group summands with the same values of $q_i = q$. The number of such summands
(for each $q$) is given by~\eqref{eq:013}. Therefore
\begin{align*}
H_{M,b} & \geq \sum_{q=1}^{M} (M+1-q) \frac{M-b-1}{M} \varphi(q) +
{\mathcal O}(M^{2+\delta})=\\
& = \frac{M-b}{M} \frac{M^3}{\pi^2} + {\mathcal
O}(M^{2+\delta})\,,
\end{align*}
where we have used~\eqref{e:euler3}. The result easily follows.
\end{proof}

%Note that the upper bound on the number of $b$-amicable pairs is
%not exact. First, $b$ must be such that the class ${\mathcal
%L}_M(\varepsilon)$ allow words with $b$ blocks $01$. This is
%solved in Lemma~\ref{l:3}. Secondly, there are exceptions among pairs $w^{(i)}$, $w^{(i+b)}$,
%$i=1,\dots,M+1-b$, which are not $b$-amicable, namely those indices $i$, such that
%the interval $J_i\cup\dots\cup J_{i+b}$ contains in its interior
%the iteration $T^{-M+1}_\varepsilon(\varepsilon)$. However, since we do not know the position of this
%point, it is difficult to determine the number of these exceptions.
%There may be other yet other pairs $w^{(i)}$, $w^{(i+b)}$, which are not $b$-amicable, namely
%if the interchange $01\leftrightarrow10$ in consecutive words
%occurs in blocks on positions that overlap. The analysis of such
%cases would however be tedious.

%%%%%%%%%%%%%%%%%%%%%%%%%%%%%%%%%%%%%%%%%%%%%%%%%%%%%%%%%%%%%%%%%%%%%%%%%%%%%%%%%%%%%%%
%\section{Asymptotic bounds on the number of 3iet factors}

%Let us now derive some asymptotic bounds on the number of factors
%occurring in any 3iet word. We use Propositions~\ref{p:zaklad}
%and~\ref{p:zaklad2} and the properties of Sturmian factors shown
%in Section~\ref{sec:sturmian}.

%%%%%%%%%%%%%%%%%%%%%%%%%%%%%%%%%%%%%%%%%%%%%%%%%%%%%%%%%%%%%%%%%%%%%%%%%%%%%%%%%%%%%%%
%\section{Lower bound on the number of 3iet factors}

The following theorem provides slightly better lower bound on $\#
3\iet(N)$ than that announced in Theorem~\ref{thm:hlavni}. Instead
of the constant $\frac13$ we give $\frac{17}{48}>\frac13$.

\begin{thm}\label{thm:dolniodhad}
The number of 3iet factors of length $N$ satisfies
$$
\# 3\iet(N)\geq \frac{17}{48\pi^2} N^4 + {\mathcal
O}(N^{3+\delta})\,,
$$
for arbitrary $\delta>0$.
\end{thm}

\begin{proof}
 Combining the formula from
 Proposition~\ref{p:zaklad2} and
Lemma~\ref{l:4}, we have
\begin{align*}
\# 3\iet(N)&\geq  2\sum_{b=0}^{\lfloor N/2\rfloor} H_{N-b,b}
\geq \frac2{\pi^2}\sum_{b=0}^{\lfloor N/2\rfloor} (N-2b)(N-b)^2 +
\sum_{b=0}^{\lfloor N/2\rfloor}{\mathcal
O}\big((N-b)^{2+\delta}\big)=\\
&= \frac2{\pi^2}\sum_{b=0}^{\lfloor
N/2\rfloor}(N^3-4N^2b+5Nb^2-2b^3)
 + {\mathcal O}(N^{3+\delta})=\\
&=  \frac{2}{\pi^2}\Big(N^3\frac{N}2 - 4N^2\frac{(N/2)^2}{2} +
 5N\frac{(N/2)^3}{3} - 2\frac{(N/2)^4}{4}\Big) + {\mathcal
 O}(N^{3+\delta}) = \\[1mm]
&=
\frac{2N^4}{\pi^2}\Big(\frac12-\frac12+\frac5{24}-\frac1{32}\Big)
+ {\mathcal O}(N^3\log N) = \frac{17}{48\pi^2} N^4 + {\mathcal
O}(N^{3+\delta})\,.\qedhere
\end{align*}
\end{proof}

\section{Precise values for small lengths}\label{last}
We have seen for Sturmian words, that given $N\in\N$, the interval
$(0,1)$ of possible values of the parameter $\varepsilon$ was
divided into $\sum_{k=1}^N\varphi(k)$ areas such that all Sturmian
words with parameter within one area had the same set of $N+1$
factors of length $N$. Similarly, one can expect that the family
of possible pairs of parameters $\varepsilon\in(0,1)$ and
$\ell\in\big(\!\max\{\varepsilon,1-\varepsilon\},1\big)$ will be
divided into regions with the same set of factors of length $N$.
%% Recall that given $\varepsilon$, the factor complexity of a 3iet
%% word $u_{\varepsilon,\ell}$ differs dependently on whether
%% $\ell\in\Z+\varepsilon\Z$ or not. Length $\ell\in\Z+\varepsilon\Z$
%% such that $\C(n)=n+\mathit{const.}$ for $n\geq N$ will certainly belong to
%% the boundaries of such regions.
The number of factors of a 3iet
word with parameters in the interior of the regions will be equal
to $2N+1$.

\begin{figure}[!ht]
\centering
\begin{minipage}{0.45\textwidth} \centering
\includegraphics{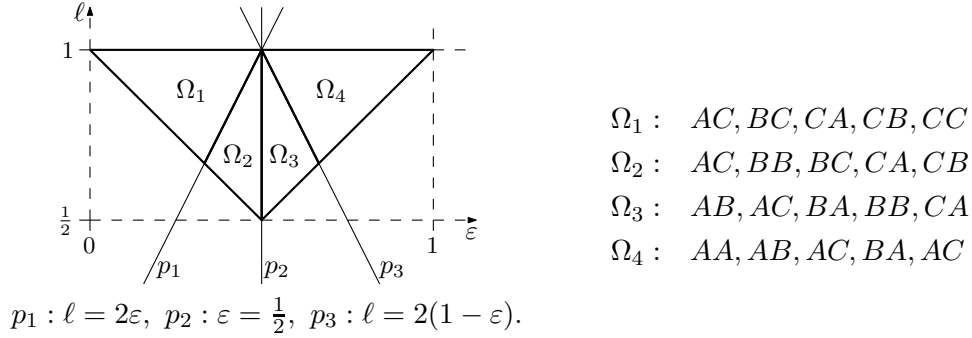}\\
$
 p_1: \ell = 2\varepsilon,\
 p_2: \varepsilon = \tfrac{1}{2},\
 p_3: \ell = 2(1-\varepsilon).
$
\end{minipage}
\begin{minipage}{0.45\textwidth}
\begin{align*}
\Omega_1: & \quad AC,BC,CA,CB,CC \\
\Omega_2: & \quad AC,BB,BC,CA,CB \\
\Omega_3: & \quad AB,AC,BA,BB,CA \\
\Omega_4: & \quad AA,AB,AC,BA,AC
\end{align*}
\end{minipage}
\caption{Division of the region of parameters $\varepsilon,\ell$ by inequalities for factors of length 2
and the corresponding lists of factors.} \label{f:3ietdelky2}
\end{figure}

Let us describe the way to obtain the list of factors of length
$N$ for fixed parameters $\varepsilon,\ell$. Denote
$d_1=\varepsilon$, $d_2=\ell-1+\varepsilon$, the discontinuity
points of the transformation $T_{\varepsilon,\ell}$. The domain of
$T_{\varepsilon,\ell}$, namely $[0,\ell)$, is divided by points
$T_{\varepsilon,\ell}^{-j}(d_1)$,
$T_{\varepsilon,\ell}^{-j}(d_2)$, $j=0,1,\dots,N-1$ generically
into $2N+1$ subintervals. (This happens for
$\ell\not\in\Z+\Z\varepsilon$ and for small $N$ also for other
$\ell$. Otherwise, some of these iterations coincide and the
number of subintervals is smaller.) Each of these subintervals
corresponds to one factor of length $N$ occurring in the language of the infinite
word $u_{\varepsilon,\ell,x_0}$ for any $x_0$. The ordering of iterations
$T_{\varepsilon,\ell}^{-j}(d_1)$, $T_{\varepsilon,\ell}^{-j}(d_2)$
in $[0,\ell)$ depends on parameters $\varepsilon,\ell$; different
orderings give rise to different lists of 3iet factors.
Values of $\varepsilon,\ell$ providing the same ordering are given
by linear inequalities in $\varepsilon,\ell$.
Figure~\ref{f:3ietdelky2} shows the division of the region of
parameters $\varepsilon,\ell$ by these inequalities for factors of length 2.
The set of
all 3iet factors of length $2$ is the union of lists given in the
figure, i.e.,
$$
3\iet(2) = \{AA,AB,AC,BA,BB,BC,CA,CB,CC\}\,,\qquad\hbox{and}\qquad
\#3\iet(2) =9\,.
$$

\begin{figure}[!ht]
\centering
\begin{minipage}{0.55\textwidth} \centering
\includegraphics{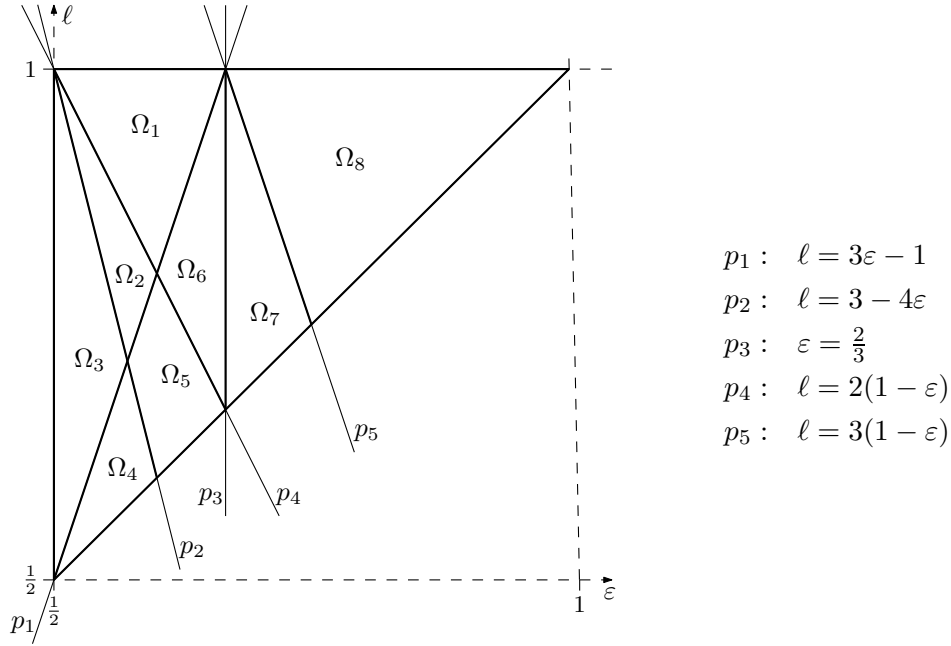}
\end{minipage}
\begin{minipage}{0.35\textwidth}
\begin{align*}
p_1: &\quad \ell = 3\varepsilon - 1 \\
p_2: &\quad \ell = 3-4\varepsilon \\
p_3: &\quad \varepsilon = \tfrac{2}{3} \\
p_4: &\quad \ell = 2(1-\varepsilon) \\
p_5: &\quad \ell = 3(1-\varepsilon)
\end{align*}
\end{minipage}
\caption{Division of the region of parameters $\varepsilon,\ell$
by inequalities for factors of length 3.} \label{f:3ietdelky3}
\end{figure}

In Figure~\ref{f:3ietdelky3}, we give the same analysis for
factors of length $N=3$. Note that the region and its division
into areas according to factors of given length must be symmetric
with respect to the axis $\varepsilon=1/2$. This corresponds to
interchange of letters $A\leftrightarrow C$ in all the factors.
Figure~\ref{f:3ietdelky3} shows only one half of the picture.

The lists of factors of length $3$ in individual areas of
Figure~\ref{f:3ietdelky3} are given as follows.
\begin{align*}
\Omega_1: & \quad AAC,ABA,ACA,BAC,CAA,CAB,CAC\\
\Omega_2: & \quad ABA,ABB,ACA,BAC,BBA,CAB,CAC\\
\Omega_3: & \quad ABB,ACA,BAC,BBA,BBB,CAB,CAC\\
\Omega_4: & \quad ABB,ACA,BAB,BAC,BBA,BBB,CAB\\
\Omega_5: & \quad ABA,ABB,ACA,BAB,BAC,BBA,CAB\\
\Omega_6: & \quad AAC,ABA,ACA,BAB,BAC,CAA,CAB\\
\Omega_7: & \quad AAB,AAC,ABA,ACA,BAA,BAB,CAA\\
\Omega_8: & \quad AAA,AAB,AAC,ABA,ACA,BAA,CAA
\end{align*}
The lists of factors in areas $\Omega_9,\dots,\Omega_{16}$ not
shown in Figure~\ref{f:3ietdelky3} are obtained from
$\Omega_1,\dots,\Omega_{8}$ by the interchange $A\leftrightarrow
C$. The set of all 3iet factors of length $3$ is therefore equal
to
\begin{align*}
3\iet(3) =\{A,B,C\}^3\setminus \{ABC,CBA\}
%\{&AAA,AAB,AAC,ABA,ABB,ACA,ACB,ACC,\\
%&BAA,BAB,BAC,BBA,BBB,BBC,BCA,BCB,BCC,\\
%&CAA,CAB,CAC,CBB,CBC,CCA,CCB,CCC\}\,,
\end{align*}
and hence $\#3\iet(3)=25$.

A computer evaluation provides the following table of values $\#
3\iet(N)$ for $N=1,\dots,10$.

\begin{center}
\begin{tabular}{c|c|c|c|c|c|c|c|c|c|c}
 $N$ & 1&2&3&4&5&6&7&8&9&10 \\ \hline
$\# 3\iet(N)$ &  3&9&25&55&113&199&339&531&809&1165 \\ \hline
\rule{0pt}{1.62em}$\displaystyle \frac{\pi^2 \# 3\iet(N)}{N^4} \approx$
&29.6&5.55&3.05&2.12&1.78&1.52&1.39&1.28&1.22&1.15
\end{tabular}
\end{center}

As we have proved above, asymptotically, the values from the last
string of the table should fall between $17/48$ and $2$. It seems
from the table that these bounds start to hold quite early, and
perhaps it is possible to improve them.

%%%%%%%%%%%%%%%%%%%%%%%%%%%%%%%%%%%%%%%%%%%%%%%%%%%%%%%%%%%%%%%%%%%%%%%%%%%%%%%%%%%%%%%
\section*{Acknowledgements}

We acknowledge financial support by the Czech Science Foundation
grant 201/09/0584 and by the grants MSM6840770039 and LC06002 of
the Ministry of Education, Youth, and Sports of the Czech
Republic. The second author was supported also by RFBR grant
09-01-00244.

%%%%%%%%%%%%%%%%%%%%%%%%%%%%%%%%%%%%%%%%%%%%%%%%%%%%%%%%%%%%%%%%%%%%%%%%%%%%%%%%%%%%%%%

%% \bibliography{../reference}

\begin{thebibliography}{10}

\bibitem{adamczewski-jtnb-14}
B.~Adamczewski.
\newblock {\em Codages de rotations et ph\'enom\`enes d'autosimilarit\'e}.
\newblock J. Th\'eor. Nombres Bordeaux {\bf 14} (2002), 351--386.

\bibitem{abmp-integers-8}
P.~Arnoux, V.~Berth{\'e}, Z.~Mas{\'a}kov{\'a}, E.~Pelantov{\'a}.
\newblock {\em Sturm numbers and substitution invariance of 3iet words}.
\newblock Integers {\bf 8} (2008), A14, 17.

\bibitem{balazi_masakova_pelantova_2}
P. Bal\'a\v zi, Z. Mas\'akov\'a, E. Pelantov\'a.
\newblock {\em Characterization of Substitution Invariant Words Coding Exchange of Three Intervals}.
\newblock Integers {\bf 8} (2008), A20, 21.

\bibitem{bdjr-dmtca}
N.~B\'edaride, E. Domenjoud, D. Jamet, J.-L. R\'emy.
\newblock {\em Number of balanced words for a given length and height}.
To appear in Discrete mathematics and theoretical computer sciences, 25 pp.

\bibitem{belov}
A. Ya. Belov, A. L. Chernyat'ev.
\newblock {\em Words with low complexity and interval exchange transformations}.
\newblock Russ. Math. Surv. {\bf 63} (2008), 158--160.

\bibitem{berstel-pocchiola-ijac-3}
J.~Berstel, M.~Pocchiola.
\newblock {\em A geometric proof of the enumeration formula for {S}turmian
  words}.
\newblock Internat. J. Algebra Comput. {\bf 3} (1993), 349--355.

\bibitem{boshernitzan_caroll}
M. D. Boshernitzan, C. R. Carroll.
\newblock {\em An extension of Lagrange's theorem to interval exchange
transformations over quadratic fields}.
\newblock J. Anal. Math. {\bf 72} (1997), 21--44.

\bibitem{cassaigne}
J. Cassaigne.
\newblock {\em Sequences with grouped factors}.
\newblock Developments in Language Theory III (DLT'97), Aristotle
University of Thessaloniki, (1998), 211--222.

\bibitem{chernyatev}
A. L. Chernyat'ev.
\newblock {\em Words with minimal growth function}.
\newblock Moscow Univ. Math. Bull. {\bf 63} (2008), 262--264.

%\bibitem{Droubay}
%X.~Droubay and G.~Pirillo.
%\newblock {\em Palindromes and sturmian words}.
%\newblock Theor. Comput. Sci. {\bf 223} (1999), 73--85.

\bibitem{ferenczi-holton-zamboni-jam-89}
S.~Ferenczi, C.~Holton, L.~Q. Zamboni.
\newblock {\em Structure of three-interval exchange transformations. {II}. {A}
  combinatorial description of the trajectories}.
\newblock J. Anal. Math. {\bf 89} (2003), 239--276.

\bibitem{gmp-jtnb-15}
L.-S. Guimond, Z.~Mas{\'a}kov{\'a}, E.~Pelantov{\'a}.
\newblock {\em Combinatorial properties of infinite words associated with
  cut-and-project sequences}.
\newblock J. Th\'eor. Nombres Bordeaux {\bf 15} (2003), 697--725.

\bibitem{hardy-wright-itn}
G.~H. Hardy, E.~M. Wright.
\newblock {\em An introduction to the theory of numbers}.
\newblock Oxford University Press, Oxford, sixth edition, (2008).
\newblock Revised by D. R. Heath-Brown and J. H. Silverman.

\bibitem{StepinKatok}
A. B. Katok, A. M. Stepin.
\newblock {\em Approximations in ergodic theory}.
\newblock Usp. Math. Nauk. {\bf 22} (1967), 81--106.

\bibitem{langevin-pmh-23}
M.~Langevin.
\newblock {\em Stimulateur cardiaque et suites de {F}arey}.
\newblock Period. Math. Hungar. {\bf 23} (1991), 75--86.

\bibitem{lipatov}
E. P. Lipatov. \newblock {\em A classification of binary collections and properties of homogeneity classes}. \newblock Problemy Kibernet. {\bf 39} (1982), 67--84 (in Russian).

\bibitem{lothaire2}
M.~Lothaire.
\newblock {\em Algebraic {C}ombinatorics on {W}ords}, volume~90 of {\em
  Encyclopedia of Mathematics and its Applications}.
\newblock Cambridge University Press, Cambridge, (2002).

\bibitem{mignosi-tcs-82}
F.~Mignosi.
\newblock {\em On the number of factors of {S}turmian words}.
\newblock Theoret. Comput. Sci. {\bf 82} (1991), 71--84.

\bibitem{morse-hedlund-ajm-62}
M.~Morse and G.~A. Hedlund.
\newblock {\em Symbolic dynamics {II}. {S}turmian trajectories}.
\newblock Amer. J. Math. {\bf 62} (1940), 1--42.

\bibitem{nathanson-gtm-195}
M.~B. Nathanson.
\newblock {\em Elementary methods in number theory}, volume 195 of {\em
  Graduate Texts in Mathematics}.
\newblock Springer-Verlag, New York, (2000).

\bibitem{fogg}
N.~Pytheas~Fogg.
\newblock {\em Substitutions in {D}ynamics, {A}rithmetics and {C}ombinatorics},
  volume 1794 of {\em Lecture Notes in Mathematics}.
\newblock Springer Verlag, (2002).

%\bibitem{vuillon}
%L.~Vuillon.
%\newblock {\em A characterization of {S}turmian words by return words}.
%\newblock European J. Combin. {\bf 22} (2001), 263--275.

%\bibitem{yasutomi}
%S. Yasutomi.
%\newblock {\em On Sturmian sequences which are invariant under some
%substitutions}.
%\newblock Number theory and its applications
%(Kyoto, 1997), Dev. Math. {\bf 2}, 347--373, Kluwer Acad. Publ.,
%Dordrecht, 1999.

\end{thebibliography}
%% \bibliographystyle{../myrefstyle}

%%%%%%%%%%%%%%%%%%%%%%%%%%%%%%%%%%%%%%%%%%%%%%%%%%%%%%%%%%%%%%%%%%%%%%%%%%%%%%%%%%%%%%%

\end{document}